\documentclass[12pt]{iopart}
\usepackage{graphicx}
\usepackage[dvips]{epsfig}
\usepackage{dcolumn}
\usepackage{subfigure}%
\usepackage{bm}
\usepackage{iopams}
\begin{document}

\title{The fractional Keller-Segel model}

\author{Carlos Escudero}

\address{Mathematical Institute, University of Oxford, 24-29 St Giles', Oxford OX1 3LB, United Kingdom}
\ead{escudero@maths.ox.ac.uk}
\begin{abstract}
The Keller-Segel model is a system of partial differential equations modelling chemotactic aggregation in cellular systems.
This model has blowing up solutions for large enough initial conditions in dimensions $d \ge 2$, but all the solutions are regular in one dimension; a mathematical fact that crucially affects the patterns that can form in the biological system. One of the strongest assumptions of the Keller-Segel model is the diffusive character of the cellular motion, known to be false in many situations. We extend this model to such situations in which the cellular dispersal is better modelled by a fractional operator. We analyze this fractional Keller-Segel model and find that all solutions are again globally bounded in time in one dimension. This fact shows the robustness of the main biological conclusions obtained from the Keller-Segel model.
\end{abstract}


\maketitle

\section{Introduction}

Biological pattern formation is a topic of growing interest in the field of applied mathematics, both due to the possibility
of developing new mathematics as well as for the broad range of important applications it might have~\cite{murray}. At the cellular level, it is known that chemotaxis plays a fundamental role in the self-organization of many biological systems.
Chemotaxis is the directed movement of an organism in response to ambient chemical gradients, that are oftenly segregated
by the cells themselves. In those cases where the chemical products are attractive (and they are therefore called chemoattractants), they lead to the phenomenon known as chemotactic aggregation: the cells accumulate in small regions of
space giving rise to high density configurations. The patterns appear therefore as strong density variations in the spatial
distribution of the cells.

One of the most important partial differential systems for understanding chemotactic aggregation is Keller-Segel model~\cite{keller}, which in a simplified form reads:
\numparts
\begin{eqnarray}
\label{keller1}
\partial_t \rho &=& D_\rho \Delta \rho - \kappa \nabla (\rho \nabla c), \\
\label{keller2}
\partial_t c &=& D_c \Delta c + \gamma \rho - \beta c.
\end{eqnarray}
\endnumparts
Here $D_\rho$ is the cellular diffusion constant, $\kappa$ the chemotactic coefficient, $\gamma$ the rate of attractant production, $\beta$ the rate of attractant depletion, $D_c$ the chemical diffusion constant, $\rho$ is the cell density, and $c$ is the chemical density. The terms
in Eq.(\ref{keller1}) include the diffusion of the cells and chemotactic drift. Eq.(\ref{keller2}) expresses the diffusion and production of attractant. This system is known to have finite time blowing up solutions for large enough initial conditions in dimensions $d \ge 2$, but all the solutions are regular for $d=1$; there is a large literature on the analysis of this problem for the Keller-Segel model as well as for some simplifications~\cite{childress,jager,nagai,herrero1,herrero2,herrero3,herrero4,osaki,hillen,horstmann}. The biological meaning of this mathematical fact was examined in Refs.~\cite{brenner,betterton}. One conclussion of these works is that in a three-dimensional system, while collapse to infinite density lines and points can occur, collapse to an infinite density sheet is mathematically impossible. Correspondingly, in a two dimensional system, it is impossible to find collapse to an infinite density line, but it is still possible to observe collapse to an infinite density point~\cite{betterton}. Let us briefly comment on this interpretation. Blow up in one dimension would mean that a real three dimensional system would be able to contract itself by using only one spatial dimension while the other two would remain invariant during this process. Global boundedness in time would invalidate this mechanism of chemotactic aggregation. Obviously, this fact crucially affects the patterns that can form. It is, however, important to point out that this connection between critical dimension and biological pattern formation is phenomenological, and a rigorous mathematical theory describing the interplay between dimensionality and geometry of the singular set for the Keller-Segel model constitutes an open problem~\cite{zaag}. So,
assuming the phenomenology, the one dimensional result is not just a model problem, but has relevance in predicting the behaviour of the biological system. Interestingly, the three-dimensional modalities of chemotactic collapse allowed by the Keller-Segel system have already been observed in experiments performed with {\it Escherichia coli}~\cite{budrene1,budrene2}. If this behaviour were general in cellular systems (collapse to points and lines, but not to sheets, in three dimensions, and to points, but not to lines, in two dimensions), and not a particular theory explaining these concrete patterns of {\it Escherichia coli}, we could name it as the \emph{Keller-Segel law} for chemotactic aggregation. This would constitute an achievement of fundamental importance in mathematical biology.

However, exceptions to this rule have already been found. In a recent {\it in vitro} experiment performed with mesenchymal cells~\cite{garfinkel}, it was shown that a two-dimensional system of these cells chemotactically aggregated into one-dimensional structures, clearly violating the Keller-Segel law. A theoretical explanation was found by taking into
account the fact that diffusion is a too unrealistic assumption for mesenchymal cells~\cite{escudero}. Actually,
these cells perform a sort of nonlocal diffusion that is better modelled by the operator $\Delta/(1-\Delta)$, defined from
its Fourier tranform
\begin{equation}
\left( \frac{\Delta}{1-\Delta} f \right) \hat{}=\frac{-k^2}{1+k^2} \hat{f},
\end{equation}
instead of the Laplacian. Performing this substitution, it was proven that the corresponding nonlocal Keller-Segel system
blows up in one dimension for large enough initial condition, leading thus to a modified Keller-Segel law that allows the type of chemotactic aggregation observed in the experiments~\cite{escudero}.

At this point there is a very important question to be answered: is the Keller-Segel law a general result or a particular property of {\it Escherichia coli} colonies? It is clear that this law is valid whenever the Keller-Segel system applies,
but what we would like to know is whether the Keller-Segel law extends beyond the validity of the partial differential system~Eqs.(\ref{keller1},\ref{keller2}). As happened in the case of mesenchymal cells, the diffusive approximation is too simplifying in many situations, particularly in those not designed at the laboratory.

Diffusion is obtained as the description of the spatiotemporal distribution of a population density of random walkers~\cite{skellam}. However, in many situations found in Nature,
L\'evy flights are commonly adopted as an efficient search strategy of living organisms~\cite{klafter,levandowsky,bartumeus}. Indeed, experimental evidence of superdiffusive behaviour has been found in biological systems~\cite{klafter,levandowsky,bartumeus}: superdiffusion is characterized by a superlinear dependence in time of the mean square displacement of the position of the dispersing population. The correct description of a population undergoing L\'evy flights is given by the substitution of the Laplacian operator by a fractional operator of the Riesz type $\Lambda^\alpha$ in one dimension, defined from its Fourier transform
\begin{equation}
(\Lambda^\alpha f)\hat{}=-|k|^\alpha \hat{f},
\end{equation}
where $1<\alpha<2$~\cite{metzler}. Direct substitution provides us in one dimension with \emph{the fractional Keller-Segel model}:
\numparts
\begin{eqnarray}
\label{fkeller1}
\partial_t \rho &=& D_\rho \Lambda^\alpha \rho - \kappa (\rho c_x)_x, \\
\label{fkeller2}
\partial_t c &=& D_c c_{xx} + \gamma \rho - \beta c.
\end{eqnarray}
\endnumparts
This model is thus presented as a more accurate description of a cellular population self-interacting chemotactically in those cases in which Brownian motion does not represent a good approximation for the population dispersal. So we can now translate our question to the present setting and ask ourselves if the Keller-Segel law is still valid for this partial differential system. It is a known mathematical fact that the Riesz operator $\Lambda^\alpha$, with $1<\alpha<2$ is less regularizing than the Laplacian, so it might be possible that the system Eqs.(\ref{fkeller1}, \ref{fkeller2}) blows up in finite time. This would imply a new breakdown of the Keller-Segel law for an important number of situations.

The goal of this article is to prove that solutions to a simplification of the system Eqs.(\ref{fkeller1}, \ref{fkeller2}) are globally bounded in time. This is an important result since it implies the robustness of the Keller-Segel
law for a large number of situations in which the original Keller-Segel system does not apply. The outline of the rest of the article is as follows. In Sec.~\ref{basic} we present the simplified model and prove the simplest estimates. In Sec.~\ref{delta} we prove the impossibility of finite time formation of Dirac masses, and in Sec.~\ref{boundedness} we prove that the solution is globally bounded in time. A summary of our results and some directions for future research are presented in Sec.~\ref{conclusions}.

\section{Nondimensionalization and basic estimates}
\label{basic}

The first step in the analysis of system Eqs.(\ref{fkeller1}, \ref{fkeller2}) is to perform a nondimensionalization. Changing variables
\numparts
\begin{eqnarray}
x &\to& \left( \frac{D_\rho}{\kappa} \right)^{1/\alpha} \hat{x}, \\
t &\to& \kappa^{-1} \hat{t}, \\
\rho &\to& \frac{\gamma D_\rho^{2/\alpha}}{D_c \kappa^{2/\alpha}} \hat{\rho},
\end{eqnarray}
\endnumparts
and supressing the hats we arrive at the nondimensional fractional Keller-Segel model:
\numparts
\begin{eqnarray}
\partial_t \rho = \Lambda^\alpha \rho - (\rho c_x)_x, \\
\frac{D_\rho^{2/\alpha}}{D_c \kappa^{2/\alpha -1}}\partial_t c = c_{xx} + \rho - \frac{D_\rho^{2/\alpha}\beta}{D_c \kappa^{2/\alpha}}c.
\end{eqnarray}
\endnumparts
Supposing that $D_c >> D_\rho$ (in the sense of much greater) and that $\kappa$ is large enough (but not making further assumptions on the value of $\beta$), we finally arrive at the system:
\numparts
\begin{eqnarray}
\label{fks1}
\partial_t \rho &=& \Lambda^\alpha \rho - (\rho c_x)_x, \\
\label{fks2}
c_{xx} &=& \delta c-\rho,
\end{eqnarray}
\endnumparts
where
\begin{equation}
\delta=\frac{D_\rho^{2/\alpha} \beta}{D_c \kappa^{2/\alpha}}.
\end{equation}
This reduction to a parabolic-elliptic system is very common and has been performed many times in the literature~\cite{jager,nagai,brenner2,velazquez}, and it can be related to direct biological measures~\cite{betterton}.

To clarify our language, let us rigorously define what we understand by chemotactic collapse:

{\sc Definition.}
{\it We say that the system (\ref{fks1}, \ref{fks2}) undergoes chemotactic collapse if there is a time $T_*<\infty$ such
that
\begin{equation}
\lim_{t \to T_*} \rho(x,t)=+\infty,
\end{equation}
for some $x \in \mathbb{R}$.}

The goal of this article is to prove that this system \emph{does not} undergo chemotactic collapse. It is important at this point to note that there are other possible ways of defining chemotactic collapse~\cite{othmer}, however, our definition is maybe more in the line of Refs.~\cite{brenner,betterton}.

We take the domain to be the full line $\mathbb{R}$, since boundary value problems for L\'evy flights are quite involved, as the long jumps make the definition of a boundary actually quite intricate (see however~\cite{zumofen,drysdale}). Suppose
that the initial condition $\rho_0$ fulfils
\begin{equation}
||\rho_0||_{L^1} < \infty,
\end{equation}
and as boundary conditions we take that $\rho$ and $c$, and their corresponding spatial derivatives, vanish when $|x| \to \infty$. Since we are dealing with densities, we are only interested in nonnegative solutions to system Eqs.(\ref{fks1},\ref{fks2}).

This system only models motion of the cells, but not any birth or death processes. This implies that the total number of cells, given by the $L^1$ norm of $\rho$, does not vary in time. It is actually easy to prove that the $L^1$ norm of $\rho$ is conserved in time:

{\sc Lemma.} {\it
Let $\rho$, $c$ be a solution to the system (\ref{fks1}, \ref{fks2}). Suppose that for the initial condition
$||\rho_0||_{L^1}<\infty$, then
\numparts
\begin{equation}
||\rho||_{L^1}(t)=||\rho||_{L^1}(0),
\end{equation}
and
\begin{equation}
||c||_{L^1}(t)=||c||_{L^1}(0).
\end{equation}
\endnumparts
}

{\it Proof.}
A direct application of the boundary conditions yields
\begin{equation}
\frac{d}{dt}||\rho||_{L^1}=\int \Lambda^\alpha \rho dx - \int (\rho c_x)_x dx=0,
\end{equation}
because $\int \Lambda^\alpha \rho dx=0$. Integrating Eq.(\ref{fks2}) we find that
\begin{equation}
||c||_{L^1}=\frac{1}{\delta}||\rho||_{L^1}.
\end{equation}
$\Box$

\section{Impossibility of $\delta$-function aggregates}
\label{delta}

The aim of this section is to prove the following

{\sc Theorem.} {\it
Let $\rho$, $c$ be a solution to the system (\ref{fks1}, \ref{fks2}). Suppose that the initial condition obeys the properties stated in the last section and also fulfils $||\rho(x,0)||_{L^2} < \infty$. Then there is a constant $C$ depending on the $L^2$ and $L^1$ norms of the initial condition such that
\begin{equation}
\max \{||\rho(\cdot,t)||_{L^2},||c(\cdot,t)||_{L^2}\} < C
\end{equation}
for all times $t > 0$.
}

{\it Proof.}
Let us start evaluating the $L^2$ norm of $\rho$:
\begin{equation}
\frac{1}{2}\frac{d}{dt}||\rho||^2_{L^2}=\int \rho \partial_t \rho dx = \int \rho \Lambda^\alpha \rho dx - \int \rho (\rho c_x)_x dx.
\end{equation}
Integration by parts yields
\begin{equation}
-\int \rho \rho_x c_x dx= \int \rho \rho_x c_x dx+ \int \rho^2 c_{xx} dx,
\end{equation}
and using Eq.(\ref{fks2}) we see that
\begin{equation}
-\int \rho \rho_x c_x dx= \frac{\delta}{2}\int \rho^2 c dx- \frac{1}{2}\int \rho^3 dx.
\end{equation}
We can use this result to find
\begin{eqnarray}
\nonumber
\frac{1}{2}\frac{d}{dt}||\rho||_{L^2}^2=\frac{1}{2}||\rho||_{L^3}^3-\frac{\delta}{2}\int \rho^2 c dx \\
-||\Lambda^{\alpha/2}\rho||_{L^2}^2 \le \frac{1}{2}||\rho||_{L^3}^3-||\Lambda^{\alpha/2}\rho||_{L^2}^2,
\label{estimater}
\end{eqnarray}
due to the positivity of $c$. The third moment of $\rho$ might be estimated as follows
\begin{equation}
||\rho||_{L^3}^3 \le ||\rho||_{L^\infty}^2||\rho||_{L^1},
\end{equation}
now choose $\chi \in (1,\alpha)$, and use a Sobolev embedding to find
\begin{equation}
||\rho||_{L^\infty}^2 \le C(||\rho||_{L^2}^2+||\Lambda^{\chi/2} \rho||_{L^2}^2).
\end{equation}
We can now invoke the Fourier transform of $\rho$
\begin{equation}
\hat{\rho}(k)=\frac{1}{\sqrt{2\pi}}\int e^{ikx}\rho(x)dx
\end{equation}
to claim that
\begin{eqnarray}
\nonumber
||\Lambda^{\chi/2}\rho||_{L^2}^2=||(\Lambda^{\chi/2}\rho)\hat{}||_{L^2}^2=\int |k|^\chi |\hat{\rho}|^2 dk= \\
\nonumber
\int_{|k|\le R}|k|^\chi |\hat{\rho}|^2 dk+\int_{|k|\ge R}\frac{|k|^\alpha}{|k|^{\alpha-\chi}}|\hat{\rho}|^2 dk \le \\
R^{\chi}||\rho||_{L^2}^2+\frac{1}{R^{\alpha-\chi}}\int |k|^\alpha |\hat{\rho}|^2 dk=R^\chi||\rho||_{L^2}^2+
\frac{1}{R^{\alpha-\chi}}||\Lambda^{\alpha/2}\rho||_{L^2}^2,
\label{chile}
\end{eqnarray}
where we have used the isometry of the Fourier transform in $L^2$. We still need to estimate the second moment of $\rho$:
\begin{eqnarray}
\nonumber
||\rho||_{L^2}^2 \le ||\rho||_{L^1}||\rho||_{L^\infty}\le \frac{1}{2\epsilon}||\rho||_{L^1}^2+ \\
\frac{\epsilon}{2}||\rho||_{L^\infty}^2\le \frac{1}{2\epsilon}||\rho||_{L^1}^2+C\frac{\epsilon}{2}
\left(||\rho||_{L^2}^2+||\Lambda^{\alpha/2}\rho||_{L^2}^2\right),
\label{smoment}
\end{eqnarray}
where we have used a Sobolev embedding. Selecting $\epsilon$ small enough we are led to conclude
\begin{equation}
||\rho||_{L^2}^2\le \left(1-C\frac{\epsilon}{2}\right)^{-1}\left(\frac{1}{2\epsilon}||\rho||_{L^1}^2+
C\frac{\epsilon}{2}||\Lambda^{\alpha/2}\rho||_{L^2}^2 \right).
\end{equation}
This inequality, in addition to Eq.(\ref{chile}) yields
\begin{eqnarray}
\nonumber
||\rho||_{L^3}^3 \le C||\rho||_{L^1}\left[(1+R^\chi)\left(1-C'\frac{\epsilon}{2}\right)^{-1} \right. \times \\
\left(\frac{1}{2\epsilon}||\rho||_{L^1}^2+
\tilde{C}\frac{\epsilon}{2}||\Lambda^{\alpha/2}\rho||_{L^2}^2 \right)
+ \left. \frac{1}{R^{\alpha-\chi}}||\Lambda^{\alpha/2}\rho||_{L^2}^2\right],
\end{eqnarray}
and substituting this result in Eq.(\ref{estimater}) we obtain
\begin{eqnarray}
\nonumber
\frac{d}{dt}||\rho||_{L^2}^2 \le C||\rho||_{L^1}\left[(1+R^\chi)\left(1-C'\frac{\epsilon}{2}\right)^{-1} \right. \times \\
\left(\frac{1}{2\epsilon}||\rho||_{L^1}^2+
\tilde{C}\frac{\epsilon}{2}||\Lambda^{\alpha/2}\rho||_{L^2}^2 \right)
+ \left. \frac{1}{R^{\alpha-\chi}}||\Lambda^{\alpha/2}\rho||_{L^2}^2\right]
-2||\Lambda^{\alpha/2}\rho||_{L^2}^2.
\end{eqnarray}
Now, by choosing a sufficiently large $R$ and a
sufficiently small $\epsilon$, we arrive at the estimate
\begin{equation}
\frac{d}{dt}||\rho||_{L^2}^2 \le C-C'||\Lambda^{\alpha/2}\rho||_{L^2}^2,
\end{equation}
for suitable constants $C,C'>0$. Repeating the steps in Eq.(\ref{smoment}) one finds
\begin{equation}
-C\frac{\epsilon_1}{2}||\Lambda^{\alpha/2}\rho||_{L^2}^2 \le
\frac{1}{2\epsilon_1}||\rho||_{L^1}^2-\left(1-C\frac{\epsilon_1}{2}\right) ||\rho||_{L^2}^2,
\end{equation}
for small enough $\epsilon_1$. Combining the last two equations yields
\begin{equation}
\frac{d}{dt}||\rho||_{L^2}^2 \le C-C'||\rho||_{L^2}^2,
\end{equation}
for suitable constants $C,C'>0$, what in turn implies the desired estimate
\begin{equation}
||\rho||_{L^2}\le C.
\end{equation}
The boundedness of the $L^2$ norm of $c$ can now be inferred from Eq.(\ref{fks2}):
\begin{equation}
\delta c=\rho+c_{xx},
\end{equation}
what implies
\begin{eqnarray}
\nonumber
\delta ||c||_{L^2}^2=\int \rho c dx+\int c c_{xx} dx \le \\
||\rho||_{L^2}||c||_{L^2}-||c_x||_{L^2}^2\le ||\rho||_{L^2}||c||_{L^2},
\end{eqnarray}
and we finally have
\begin{equation}
||c||_{L^2} \le \frac{||\rho||_{L^2}}{\delta}.
\end{equation}
$\Box$

{\sc Corollary.} {\it
This result precludes the formation of Dirac masses in the system (\ref{fks1}, \ref{fks2}).
}

\section{Global boundedness in time}
\label{boundedness}

We will prove in this section the global boundedness in time of the $L^\infty$ norm of $\rho$ and $c$:

{\sc Theorem.} {\it
Let $\rho$, $c$ be a solution to the system (\ref{fks1}, \ref{fks2}). Suppose that the initial condition obeys the properties stated in the last two sections and also fulfils
\begin{equation}
\left|\left|\frac{d \rho(x,0)}{dx}\right|\right|_{L^2} < \infty.
\end{equation}
Then there is a constant $C$ depending on the $L^2$ norm of the first spatial derivative of the
initial condition, and the $L^1$ and $L^2$ norms of the initial condition such that
\begin{equation}
\max_x \{\rho(x,t),c(x,t)\} < C
\end{equation}
for all times $t > 0$.
}

{\it Proof.}
Let us start estimating the time evolution of the $L^2$ norm of $\rho_x$:
\begin{eqnarray}
\nonumber
\frac{1}{2}\frac{d}{dt}||\rho_x||_{L^2}^2=\int \rho_x \partial_t \rho_x dx=\int \rho_x \Lambda^\alpha \rho_x dx- \\
\int \rho_x(\rho c_x)_{xx} dx= -\int \rho_x(\rho c_x)_{xx} dx-||\Lambda^{\alpha/2}\rho_x||_{L^2}^2.
\label{estimaterx}
\end{eqnarray}
The last integral in this equation may be rearranged in the following way
\begin{eqnarray}
\nonumber
I=-\int \rho_x(\rho c_x)_{xx}dx=-\int \rho_x \rho_{xx}c_x dx \\
-2\int \rho_x^2 c_{xx} dx-\int \rho_x \rho c_{xxx} dx=I_1+I_2+I_3,
\end{eqnarray}
and we can manipulate each of the terms by reiteratively using Leibniz's rule, integration by parts, and the boundary
conditions:
\begin{equation}
I_1=\frac{\delta}{2}\int \rho_x^2c dx-\frac{1}{2}\int \rho_x^2\rho dx,
\end{equation}
\begin{equation}
I_2=-2\delta \int \rho_x^2 c dx+ 2\int \rho_x^2 \rho dx,
\end{equation}
and
\begin{equation}
I_3=\frac{\delta^2}{2}\int \rho^2 c dx-\frac{\delta}{2}\int \rho^3 dx+ \int \rho \rho_x^2 dx.
\end{equation}
We finally find
\begin{eqnarray}
\nonumber
I=I_1+I_2+I_3=\frac{5}{2}\int \rho_x^2 \rho dx+ \frac{\delta^2}{2}\int \rho^2 c dx- \\
\frac{3\delta}{2}\int \rho_x^2 c dx
-\frac{\delta}{2}\int \rho^3 dx \le \frac{5}{2}\int \rho_x^2 \rho dx+ \frac{\delta^2}{2}\int \rho^2 c dx,
\end{eqnarray}
due to the positivity of $\rho$ and $c$. Substituting in Eq.(\ref{estimaterx}) we get
\begin{equation}
\frac{1}{2}\frac{d}{dt}||\rho_x||_{L^2}^2 \le \frac{5}{2}\int \rho_x^2 \rho dx+ \frac{\delta^2}{2}\int \rho^2 c dx
-||\Lambda^{\alpha/2}\rho_x||_{L^2}^2.
\end{equation}
It follows immediately that the second integral can be simply estimated by
\begin{equation}
\int \rho^2 c dx \le ||\rho||_{L^2}^2||c||_{L^\infty},
\end{equation}
while the first integral might be estimated as follows
\begin{equation}
\int \rho_x^2 \rho dx \le ||\rho_x||_{L^\infty}||\rho_x||_{L^2}||\rho||_{L^2}
\le \frac{1}{2\epsilon}||\rho_x||_{L^2}^2||\rho||_{L^2}^2+2\epsilon ||\rho_x||_{L^\infty}^2.
\end{equation}
We can use the Sobolev embedding
\begin{equation}
||\rho_x||_{L^\infty}^2 \le C\left(||\rho_x||_{L^2}^2+||\Lambda^{\alpha/2}\rho_x||_{L^2}^2\right),
\end{equation}
to arrive at
\begin{eqnarray}
\label{ile}
\nonumber
I \le \frac{\delta^2}{2}||\rho||_{L^2}^2||c||_{L^\infty}+ \\
\frac{5}{2}\left[\frac{1}{2\epsilon}||\rho_x||_{L^2}^2||\rho||_{L^2}^2
+2\epsilon C\left(||\rho_x||_{L^2}^2+ ||\Lambda^{\alpha/2}\rho_x||_{L^2}^2\right) \right],
\end{eqnarray}
so Eq.(\ref{estimaterx}) finally reads
\begin{eqnarray}
\nonumber
\frac{1}{2}\frac{d}{dt}||\rho_x||_{L^2}^2 \le
\frac{\delta^2}{2}||\rho||_{L^2}^2||c||_{L^\infty}
+\frac{5}{2}\left[\frac{1}{2\epsilon}||\rho_x||_{L^2}^2||\rho||_{L^2}^2 \right. \\
\left. +2\epsilon C\left(||\rho_x||_{L^2}^2+ ||\Lambda^{\alpha/2}\rho_x||_{L^2}^2\right) \right]
-||\Lambda^{\alpha/2}\rho_x||_{L^2}^2.
\label{estimaterx2}
\end{eqnarray}
We still need to estimate the $L^\infty$ norm of $c$. To this end we will use the Fourier transformed version of
Eq.(\ref{fks2}):
\begin{equation}
-k^2\hat{c}=\delta \hat{c}-\hat{\rho},
\end{equation}
which implies
\begin{equation}
|\hat{c}|^2=\frac{|\hat{\rho}|^2}{(k^2+\delta)^2}.
\end{equation}
This allows us to estimate the $L^2$ norm of $c_x$
\begin{eqnarray}
\nonumber
||c_x||_{L^2}^2=\int k^2|\hat{c}|^2dk=\int \frac{k^2}{(k^2+\delta)^2}|\hat{\rho}|^2dk \le \\
||\rho||_{L^1}^2
\int \frac{k^2}{(k^2+\delta)^2}dk \le C.
\end{eqnarray}
This fact, in addition to the Sobolev embedding
\begin{equation}
||c||_{L^\infty}^2 \le C\left(||c||_{L^2}^2+||c_x||_{L^2}^2 \right),
\end{equation}
yields the desired inequality
\begin{equation}
||c||_{L^\infty} \le C.
\end{equation}
This result, together with Eq.(\ref{estimaterx2}), yields, for $\epsilon$ small enough and suitable
positive constants $C_1$, $C_2$, and $C_3$, the inequality
\begin{equation}
\label{estimaterxf}
\frac{d}{dt}||\rho_x||_{L^2}^2 \le C_1 + C_2 ||\rho_x||_{L^2}^2 -C_3 ||\Lambda^{\alpha/2}\rho_x||_{L^2}^2.
\end{equation}
To finish our proof we still need the estimate
\begin{eqnarray}
\nonumber
||\rho_x||_{L^2}^2=\int |k|^2 \hat{\rho}^2 dk = \int_{|k|<R} |k|^2 \hat{\rho}^2 dk +
\int_{|k|>R} \frac{|k|^{2+\alpha}}{|k|^\alpha} \hat{\rho}^2 dk \le \\
R^2||\rho||_{L^2}^2 + \frac{1}{R^\alpha}||\Lambda^{\alpha/2}\rho_x||_{L^2}^2,
\end{eqnarray}
or, what is the same,
\begin{equation}
-||\Lambda^{\alpha/2}\rho_x||_{L^2}^2 \le R^{\alpha+2} ||\rho||_{L^2}^2 - R^{\alpha} ||\rho_x||_{L^2}^2.
\end{equation}
By choosing $R$ large enough, and combining this last inequality together with Eq.(\ref{estimaterxf}) we obtain
\begin{equation}
\frac{d}{dt}||\rho_x||_{L^2}^2 \le C_1 - C_2 ||\rho_x||_{L^2}^2,
\end{equation}
for suitable positive constants $C_1$ and $C_2$, what implies
\begin{equation}
||\rho_x||_{L^2}^2 \le C.
\end{equation}
We recall again the Sobolev embedding
\begin{equation}
||\rho||_{L^\infty}^2 \le C\left(||\rho||_{L^2}^2+||\rho_x||_{L^2}^2\right),
\end{equation}
to conclude
\begin{equation}
||\rho||_{L^\infty} \le C.
\end{equation}
$\Box$

{\sc Corollary.} {\it
This result prohibits the chemotactic collapse for the system Eqs. (\ref{fks1}, \ref{fks2}).
}

\section{Conclusions}
\label{conclusions}

In this work we have introduced and analyzed the fractional Keller-Segel model. This model is thought to be an extension of the standard Keller-Segel model to such situations in which the movements of the cells cannot be described by random walks, and a L\'evy flight description fits much better the cellular trajectories. Our analysis has revealed that solutions to the fractional Keller-Segel model exist globally in time in one dimension, as it happens with the standard Keller-Segel model, implying that chemotactic collapse to a sheet cannot occur in a three-dimensional setting, and chemotactic collapse to a line is impossible in two dimensions. This fact shows the robustness of the Keller-Segel law, that extends beyond its initial range of validity and applies to a much broader class of situations.

We have performed our analysis taking the domain of the fractional Keller-Segel model to be the whole line $\mathbb{R}$,
which is a good first approximation (in the standard Keller-Segel model the proof of existence do not vary much from the unbounded domain to the bounded one~\cite{childress}). However, it would be interesting to analize our model in a bounded domain, specially with Neumann boundary conditions, due to its biological relevance in experiments~\cite{childress,escudero,velazquez}. As shown in the diffusive case, homogeneous solutions become unstable for certain range of the parameter values, and a rich bifurcation structure, which may be supercritical or subcritical, arises~\cite{childress}.

Finally, we would like to remark that all the equations that have been written in this work are a continuous description of a discrete process. The number of cells involved must be integer, a fact that becomes crucially important when this number is not high enough. This type of discrete processes admits a continuous description by means of stochastic partial differential equations~\cite{escudero2}, that favours its analytical assessment. It would be interesting to study the effect of such stochastic corrections on the problem of chemotactic collapse, to see whether or not they are able to modify the critical dimension of the problem. These and other considerations will be the object of future research.

\section*{Acknowledgments}

This work has been partially supported by the Ministerio de Educaci\'on y Ciencia (Spain) through Projects No. EX2005-0976 and FIS2005-01729.

\end{document}